\documentclass[12pt]{article}
\usepackage{amssymb,amsmath,latexsym}
\usepackage{amsthm}
\usepackage{fontenc}
\numberwithin{equation}{section}
\begin{document}
\author{Alexander E Patkowski}
\title{Remarks on a Bailey pair with one free parameter}
\date{\vspace{-5ex}}
\maketitle
\abstract{We offer a more general Bailey pair than one that was proved in two different papers by two different methods [5, 12].}

\section{Introduction}
Recall that a Bailey pair $(\alpha_n(a,q), \beta_n(a,q))=(\alpha_n, \beta_n)$ is a pair of sequences that satisfy (relative to $a$) [6]
\begin{equation}\beta_n=\sum_{i\ge0}\frac{\alpha_i}{(q)_{n-i}(aq)_{n+i}}.\end{equation}
(Refer to [10] for $q$-series notation.) In [22, Lemma 2.2] we offered a new Bailey pair for 
\begin{equation}\beta_n=\frac{(-1)^nq^{n^2}}{(-a)_{2n}(q^2;q^2)_n},\end{equation}
where $a=1,$ or $a=q.$ As a direct consequence of this pair, we were able to offer new information on the distinct rank parity
function \begin{equation}\sigma(q)=\sum_{n\ge0}\frac{q^{n(n+1)/2}}{(1+q)(1+q^2)\cdots(1+q^n)}.\end{equation}
The function $\sigma(q)$ was shown in [1] to be related to the arithmetic of $\mathbb{Q}(\sqrt{6}),$ and therefore $\mathbb{Q}(\sqrt{2}),$
and $\mathbb{Q}(\sqrt{3}).$ The key to this observation was the use of Bailey pairs to relate $\sigma(q)$ to indefinite quadratic forms.
\par Andrews [5, pg.72, eq.(45)] also offered a proof of the pair with the $\beta_n$ given in (1.2) with $a=q$ using a recurrence approach, and subsequently establishing new 
partition theorems using different limiting cases of Bailey's lemma than used in [12].

\section{The Bailey pair} 
We will apply the same proof offered in [12] but in greater generality.
\\*
{\bf Theorem 2.1} \it With free parameter $x,$ we have the Bailey pair $(\alpha_n'(a^2,x,q^2), \beta_n'(a^2,x,q^2))$ where
\begin{equation}\alpha_n'(a^2,x,q^2)=\frac{q^{n^2-n}(-x)^n(1-a^2q^{4n})(a^2/x;q^2)_n}{(1-a^2)(q^2x;q^2)_n}\sum_{0\le j\le n}\frac{(1-aq^{2j-1})(a)_{j-1}(x;q^2)_j}{q^{j(j-1)/2}x^j(q)_j(a^2/x;q^2)_j},\end{equation}
\begin{equation}\beta_n'(a^2,x,q^2)=\frac{(-1)^nq^{n^2}(1-x)}{(-a)_{2n}(q^2;q^2)_n(1-xq^{2n})}.\end{equation}
\rm.
\begin{proof}
From [2, Theorem 2.3] we have the pair $(\alpha_n, \beta_n)$ relative to $a$ where

\begin{equation}\alpha_n(a,b,c,q)=\frac{q^{n^2}(bc)^n(1-aq^{2n})(a/b)_n(a/c)_n}{(1-a)(bq,cq)_n}\sum_{0\le j\le n}\frac{(-1)^j(1-aq^{2j-1})(a)_{j-1}(b,c)_j}{q^{j(j-1)/2}(bc)^j(q,a/b,a/c)_j},\end{equation}
\begin{equation}\beta_n(a,b,c,q)=\frac{1}{(bq,cq)_n}.\end{equation}
\rm
Recall from [12, eq.(2.10)--(2.11)] the result,
\begin{equation}\bar{\alpha}_n(a^2,q^2)=\frac{(1+aq^{2n})}{(1+a)q^n}\alpha_n(a,q),\end{equation}
\begin{equation}\bar{\beta}_n(a^2,q^2)=\frac{q^{-n}}{(-a;q)_{2n}}\sum_{k\ge0}\frac{(-1)^{n-k}q^{(n-k)^2-(n-k)}}{(q^2;q^2)_{n-k}}\beta_k(a,q).\end{equation}
The equations (2.5)--(2.6) allow us to change the base of a Bailey pair from $q$ to $q^2.$
\par
From Fine [9, pg.18, eq.(16.3) with $a\rightarrow0$] we have
\begin{equation}\sum_{i\ge0}\frac{(-1)^{n-i}q^{\{(n-i)^2-(n-i)\}/2}}{(q)_{n-i}(bq)_i}=\frac{(-1)^nq^{n(n+1)/2}(1-b)}{(q)_n(1-bq^n)}.\end{equation}
Now putting $b=-c$ in (2.3)--(2.4), replacing $c$ by $\sqrt{x},$ and then inserting the resulting pair in (2.5)--(2.6) gives Theorem 2.1 after noting (2.7).
\end{proof}
\par
We now offer some corollaries as special cases of Theorem 2.1 that will be noted in the next section.
\\*
{\bf Corollary 2.2} \it We have the Bailey pair 
\begin{equation}\alpha_n'(q^2,-q,q^2)=\frac{(-1)^nq^{n(n-1)/2}(1-q^{2n+1})}{(1-q)},\end{equation}
\begin{equation}\beta_n'(q^2,-q,q^2)=\frac{(-1)^nq^{n^2}}{(-q)_{2n+1}(q^2;q^2)_n}.\end{equation}
\\*
{\bf Corollary 2.3} \it We have the Bailey pair 
\begin{equation}\alpha_n'(q^4,q,q^2)=\frac{(-1)^nq^{n^2}(1-q^{4n+4})}{(1-q^4)}\sum_{0\le j\le n}q^{-j(j+1)/2},\end{equation}
\begin{equation}\beta_n'(q^4,q,q^2)=\frac{(-1)^nq^{n^2}}{(-q)_{2n}(q^2;q^2)_n(1-q^{4n+2})}.\end{equation}
\\*
{\bf Corollary 2.4} \it We have the Bailey pair 
\begin{equation}\alpha_n'(q^2,q,q^2)=\frac{(-1)^nq^{n^2}(1+q^{2n+1})}{(1-q^2)}\left(\sum_{0\le j\le n}q^{-j(j+1)/2}+\sum_{0\le j\le n-1}q^{-j(j+1)/2}\right),\end{equation}
\begin{equation}\beta_n'(q^2,q,q^2)=\frac{(-1)^nq^{n^2}}{(-q)_{2n}(q^2;q^2)_n(1-q^{2n+1})}.\end{equation}

\section{Some partitions and $q$-series} \rm
We will use some special instances of Bailey's lemma [6]
\begin{equation}\sum_{n\ge0}(X)_n(Y)_n(aq/XY)^n\beta_n=\frac{(aq/X)_{\infty}(aq/Y)_{\infty}}{(aq)_{\infty}(aq/XY)_{\infty}}\sum_{n\ge0}\frac{(X)_n(Y)_n(aq/XY)^n\alpha_n}{(aq/X)_n(aq/Y)_n}.\end{equation}

In [1] we find Andrews et al. considered the function
\begin{equation} \sigma^{*}(q)=\sum_{n\ge1}\frac{(-1)^nq^{n^2}}{(q;q^2)_n},\end{equation}
which generates $O(n),$ the number of partitions of $n$ into odd parts with the property that if a number appears then all smaller numbers appear as parts as well, weighted
by $-1$ if the largest part is congruent to $1\pmod{4}$ and $+1$ if the largest part is congruent to $3\pmod{4}.$ It was noted in [1] that $\sigma^{*}(q)$ is also related to $\mathbb{Q}(\sqrt{6}).$ Further notes on $\sigma^{*}(q)$ can be found in [11], and more examples related to real quadratic forms are given in [5, 7, 8, 11, 12]. We consider a similar function
\begin{equation}\sum_{n\ge1}\frac{(-1)^nq^{n^2}}{(q;q^2)_n(1+q^{2n-1})}=\sum_{n\ge1}\frac{(-1)^nq^{1+3+\cdots+2n-1}}{(1-q)(1-q^3)\cdots(1-q^{2n-3})(1-q^{(2n-1)+(2n-1)})}.\end{equation}
The $q$-series in (3.3) generates the same partitions counted by $O(n)$ (and same weight function) with the additional condition that the largest part appears an odd number of times. We denote such partitions to be $O^{*}(n).$ \\*
{\bf Corollary 3.5} We have
\begin{equation}2\sum_{n\ge1}O^{*}(n)(-q)^n=\sum_{n\ge1}q^{n^2}\sum_{j=-n}^{n-1}q^{-j(j+1)/2},\end{equation}
\begin{proof} Take the Bailey pair in Corollary 2.3 and insert it into (3.1) (with $a=q^2,$ then $q\rightarrow q^2$) with $X=q^2,$ $Y=-q^2.$ The proof is complete
after multiplying both sides by $2$ upon noting that $2\sum_{0\le j\le n}q^{-j(j+1)/2}=\sum_{j=-n-1}^{n}q^{-j(j+1)/2}$ and then shift the resulting indefinite quadratic sum over $n$ with $n\rightarrow n-1.$\end{proof}
Set $L=\mathbb{Q}(\sqrt{2}),$ let $O_L$ be the ring of integers of $L,$ and let $\breve{a}\subset O_L$ denote an ideal. For such an ideal we denote its norm function to be
$N(\breve{a}).$ Corollary 3.5 may be used to relate $O^{*}(n)$ to the number of inequivalent elements of $O_L$ with norm $N(\breve{a})=2x^2-y^2.$ \\*
{\bf Corollary 3.6} \it We have
\begin{equation}2q^{-1}\sum_{n\ge1}O^{*}(n)q^{8n}=\sum_{\substack{\breve{a}\subset O_L \\ N(\breve{a})\equiv-1\pmod{8}}}(-1)^{\frac{N(\breve{a})+1}{8}}q^{N(\breve{a})},\end{equation} \rm
\begin{proof} We use [1, Lemma 3] and Corollary 3.5. Since $2(-1)^nO^{*}(n)$ is equal to the number of solutions of $n=i^2-j(j+1)/2$ with $-i\le j \le i-1,$ $i\ge1,$ we may write $2(2i)^2-(2j+1)^2=8n-1.$ Any solution of 
$N(\breve{a})=2x^2-y^2=8n-1,$ $n\in\mathbb{N},$ must have $y$ odd, and subsequently $x$ even. Further, if we write $x=2i,$ $y=2j+1,$ we have
\begin{equation} -(2i)< (2j+1)\le (2i), i>0.\end{equation}
The  solutions of (3.6) are precisely the pairs $(j,i)$ that satisfy $-i\le j \le i-1,$ $i\ge1.$ \end{proof} 
{\bf Corollary 3.7 [3, Entry 9.4.3.]} \it We have
\begin{equation}\sum_{n\ge0}\frac{q^{n(2n+1)}}{(-q)_{2n+1}}=\sum_{n\ge0}q^{n(3n+1)/2}(1-q^{2n+1}).\end{equation}
\rm 
\begin{proof}
Insert the Bailey pair in Corollary 2.2 into (3.1) with $X=q^2,$ $Y\rightarrow\infty.$\end{proof}
We mention this result in passing only to emphasize that Theorem 2.1 contains a broad range of identities. Corollary 3.7 has appeared in [3, pg.233, eq.(9.4.4)], and was noted in [12] due to its
relevance to the $q$-series
$$\sum_{n\ge0}\frac{q^{n(2n+1)}}{(-q)_{2n}},$$ which was found to be lacunary and related to $\sigma(q)$ in [12], by using the $x\rightarrow0$ instance of Theorem 2.1.
\par In [7] Bringmann and Kane consider the $q$-series 
\begin{equation}f_1(q)=\sum_{n\ge0}\frac{q^{n(n+1)/2}}{(-q)_n(1-q^{2n+1})},\end{equation}
and related it to the arithmetic of $\mathbb{Q}(\sqrt{2}).$ We offer some further information on $f_1(q).$ \\*
{\bf Corollary 3.8} \it We have,
\begin{equation}\sum_{n\ge0}\frac{q^{n(n+1)}}{(-q^2;q^2)_n(1-q^{2n+1})}=\sum_{n\ge0}q^{n(n+1)}(1+q^{2n+2})\sum_{j=0}^{n}q^{-j(j+1)/2}.\end{equation}
Therefore, the ``even function" of $f_1(q^2)$ is lacunary.
\begin{proof} Let $f_1'(q)$ be the left side of (3.9). Then clearly, $(f_1'(q)+f_1'(-q))/2=f_1(q^2).$ The result follows after inserting the Bailey pair
in Corollary 2.3 into the $X=q^2,$ $Y=-q^3$ instance of (3.1) to get (3.9). \end{proof} \rm
Another proof may be obtained using Corollary 2.4, and we leave this to the reader. \\*
\rm
{\bf Corollary 3.9} \it Let $s_{+}t(n)$ be the number of representations of $n$ as a sum of a triangular number $i$ and a square $j$ weighted by $(-1)^j,$ or
$$s_{+}t(n)=\sum_{\substack{r\in\mathbb{Z},k\ge0 \\ n=r^2+k(k+1)/2}}(-1)^{r^2}.$$
 Then,
\begin{equation}\sum_{n\in\mathbb{Z}, m\ge0}(-1)^nq^{n^2+m(m+1)/2}=\sum_{n\ge0}s_{+}t(n)q^n=\sum_{n\ge0}(-1)^nq^{n(n+1)}\left(\sum_{0\le j\le n}q^{-j(j+1)/2}+\sum_{0\le j\le n-1}q^{-j(j+1)/2}\right).\end{equation}
\rm \begin{proof} We use the $X=-q,$ $Y=-q^2$ case of (3.1) coupled with the Bailey pair in Corollary 2.4 and invoke the identity (a special limiting case of [9, pg.18, eq.(16.3)])
\begin{equation}\sum_{n\ge0}\frac{(-1)^nq^{n(n+1)}(1-x)}{(q^2;q^2)_n(1-xq^{2n})}=\frac{(q^2;q^2)_{\infty}}{(xq^2;q^2)_{\infty}}. \end{equation} \end{proof}
We note in closing that Corollary 3.9 gives a mapping between the number of inequivalent elements of $O_L$ with norm $8n+1$ to 
the number of inequivalent elements of $O_{L'},$ where $L'=\mathbb{Q}(\sqrt{-2}),$ with norm $8n+1.$ \par We also mention it is possible to obtain a more general expansion for the 
product \begin{equation}\frac{(q)_{\infty}(q^2;q^2)_{\infty}}{(-q)_{\infty}(xq^2;q^2)_{\infty}}=\sum_{n\ge0}q^{n^2}(-x)^n(1-q^{2n+1})\frac{(q^2/x;q^2)_n}{(q^2x;q^2)_n}\sum_{0\le j\le n}\frac{(1+q^{j})(x;q^2)_j}{q^{j(j-1)/2}x^j(q^2/x;q^2)_j},\end{equation} using the same limiting case of (3.1) with the $a=q$ case of Theorem 2.1. A nice corollary of (3.12) with $x\rightarrow0$ is the famous expansion
due to Rogers [13] for the weight $1$ modular form
\begin{equation}\prod_{n\ge1}(1-q^n)^2=\sum_{n\ge0}q^{n(2n+1)}(1-q^{2n+1})\sum_{|j|\le n}(-1)^jq^{-j(3j+1)/2}.\end{equation}

1390 Bumps River Rd. \\*
Centerville, MA
02632 \\*
USA \\*
E-mail: alexpatk@hotmail.com
\end{document}